\documentclass[12pt]{article}
\usepackage{amsmath,amsfonts}
\usepackage{graphicx}
\usepackage{float}
\usepackage{subcaption}
\usepackage{appendix}
\usepackage{color}
\usepackage{suffix}

\newcommand {\eq} [1] {\begin{equation}\label{#1}}
\newcommand {\en} {\end{equation}}
\newcommand {\cB}       {{\mathcal B}}

\newcommand {\cS}       {{\mathcal S}}
\newcommand {\cT}       {{\mathcal T}}

\newcommand {\proof} {\par{\it Proof}. \ignorespaces}
\newcommand {\eproof}
      {\space
        {\ \vbox{\hrule\hbox{\vrule height1.3ex\hskip0.8ex\vrule}\hrule}}
        \par}
\newcommand {\mat}      [1] {\left[\begin{array}{#1}}
\newcommand {\rix}          {\end{array}\right]}
\newtheorem{theorem}           {Theorem}
\newtheorem{lemma}    [theorem]{Lemma}
\newtheorem{definition}    [theorem]{Definition}
\newtheorem{corollary}[theorem]{Corollary}

\newtheorem{remark}            {Remark}

\newcommand {\diag}     {\mathop{\rm diag}\nolimits}
\newcommand {\rank}     {\mathop{\rm rank}\nolimits}

 \catcode`@=11
 \font\tenex=cmex10 
 \newdimen\p@renwd
 \setbox0=\hbox{\tenex B} \p@renwd=\wd0 
 \def\bmat#1{\begingroup \m@th
   \setbox\z@\vbox{\def\cr{\crcr\noalign{\kern2\p@\global\let\cr\endline}}%
     \ialign{$##$\hfil\kern2\p@\kern\p@renwd&\thinspace\hfil$##$\hfil
       &&\quad\hfil$##$\hfil\crcr
       \omit\strut\hfil\crcr\noalign{\kern-\baselineskip}%
       #1\crcr\omit\strut\cr}}%
   \setbox\tw@\vbox{\unvcopy\z@\global\setbox\@ne\lastbox}%
   \setbox\tw@\hbox{\unhbox\@ne\unskip\global\setbox\@ne\lastbox}%
   \setbox\tw@\hbox{$\kern\wd\@ne\kern-\p@renwd\left[\kern-\wd\@ne
     \global\setbox\@ne\vbox{\box\@ne\kern2\p@}%
     \vcenter{\kern-\ht\@ne\unvbox\z@\kern-\baselineskip}\,\right]$}%
   \null\;\vbox{\kern\ht\@ne\box\tw@}\endgroup}
 \catcode`@=12
\begin{document}

\title{Stabilization of linear Port-Hamiltonian Descriptor Systems via Output Feedback}
\author{Delin Chu\footnotemark[1]
\and Volker Mehrmann\footnotemark[2]}

\maketitle

\renewcommand{\thefootnote}{\fnsymbol{footnote}}
\footnotetext[1]{Department of Mathematics, National University of Singapore, Singapore 119076. Email: {\tt
 matchudl@nus.edu.sg}.
 {\color{black} D. C. was supported by  NUS Research Grant A-8000430-00-00 and a Humboldt Fellowship at TU Berlin (Germany)}     }
\footnotetext[2]{Institute of Mathematics, TU Berlin, Stra{\ss}e des 17.~Juni 136, 10623 Berlin, Germany.  Email: {\tt mehrmann@math.tu-berlin.de}.
{\color{black} V. M. has been supported by Deutsche Forschungsgemeinschaft (DFG) through
 the SPP1984 “Hybrid and multimodal energy systems” Project: Distributed
 Dynamic Security Control and by the DFG Research Center Math+, Project
 No. 390685689. Advanced Modeling, Simulation, and Optimization of Large
 Scale Multi-Energy Systems}            }

\begin{abstract} {\color{black} The  structure preserving stabilization  of (possibly non-regular) linear port-Ha\-mil\-to\-nian descriptor (pHDAE) systems
by  output feedback is discussed. For general descriptor systems the characterization when there exist output feedbacks that lead to an asymptotically stable closed loop system is a very hard and partially an open problem.} In contrast to this it is shown that for systems in pHDAE representation this problem can be completely solved. Necessary and sufficient conditions are presented that guarantee that there exist a  proportional {\color{black} and/or derivative} output feedback such that the resulting closed-loop port-Hamiltonian descriptor system is  asymptotically stable. For this it is also necessary that the output feedback also makes the problem regular and of index at most one.  A complete characterization when this is possible is  presented as well.
\end{abstract}



\section{Introduction}\label{intro}
In this paper we study output feedback controls to make a descriptor system, often called differential-algebraic system (DAE) asymptotically stable.  Consider a \emph{general descriptor system} of the form
\begin{eqnarray}
E \dot x &=&A x + Bu, \ x(t_0)=x_0\nonumber\\
y&=& C x  + D u, \label{gendae}
\end{eqnarray}
with $E,A\in {\mathbb C}^{\ell,n}$,
$B\in{\mathbb C}^{\ell,m}$, $C\in{\mathbb C}^{p,n}$ $D\in {\mathbb C}^{m,m}$.  Here ${\mathbb C}^{p,n}$ denotes the complex $p\times n $ matrices, $u$ is the input, $y$ is the output and $x$ is the generalized state (descriptor) vector,
and $\dot x$ denotes the time derivative.

We formulate our results for complex systems but the results hold analogously for systems with real coefficients. In the following, the real part of a complex number $z$ is denote by $\Re(z)$ and we denote that a Hermitian matrix $M$ is positive semidefinite (positive definite) by $M\geq 0$ ($M>0$).

In our analysis and in the construction of feedbacks, we need to {\color{black} perform equivalence transformations for the system}. For general descriptor systems, these are changes of bases $x=T\tilde x$, $u=V \tilde u$, $y= Y \tilde y$ and multiplications of the state equation by $S$, where the matrices $S,T,V,Y$ are invertible.

{\color{black} The spectral properties of the matrix pencil $\lambda E-A$ associated with general descriptor systems of the form \eqref{gendae}
are characterized via the  Kronecker canonical form \cite{Gan59b}}.
A value $\lambda_0\in\mathbb C$ is called a (finite) eigenvalue of $\lambda E-A$ if $
\operatorname{rank}(\lambda_0E-A)<\max_{\alpha\in\mathbb C}
\operatorname{rank}(\alpha E-A)$.
Furthermore,
$\lambda_0=\infty$ is said to be an eigenvalue of $\lambda E-A$ if zero is an eigenvalue of $\lambda A-E$.
The size of the largest Jordan block {\color{black} associated with the eigenvalue $\infty$}
is called the \emph{index} $\nu$ of the pencil $\lambda E-A$, where, by convention,  $\nu=0$ if $E$ is invertible.
The matrix pencil $\lambda E-A$ is called \emph{regular} if $\ell=n$ and
$\operatorname{det}(\lambda_0 E-A)\neq 0$ for some $\lambda_0 \in \mathbb C$,
otherwise it is called \emph{singular}.
For a given input $u$, an initial condition $x_0$ is called \emph{consistent} if the initial value problem has at least one classical solution.

When descriptor systems are generated in an automated modularized modeling framework such as e.g. \cite{Fri03}, then the resulting system typically is an over- or underdetermined (singular) system. For such singular systems, existence and uniqueness of the solutions for a given control input and  given consistent initial values $x(t_0)=x_0$ {\color{black} can only be guaranteed}  if $E,A$ are square and the pencil $\lambda E-A$ is \emph{regular}. If this is not the case then a regularization or reformulation is necessary, see \cite{CamKM12,KunM06}. In control design this is often  done via state or output feedback, see e.g. \cite{BunBMN99,ChuMN99}. Feedback design is also used classically to make the system  asymptotically stable \cite{Kai80,ZhoDG96}.
However, to do this with output feedback is a difficult and partially open problem even if $E=I$, the identity matrix, see e.g.~\cite{BriGSS21,SyrADG97}.

Note that for {\color{black} descriptor systems} the definition of stability and asymptotic stability is not defined in a uniform way in the literature.
Some authors just require that the finite eigenvalues of $\lambda E-A$ are in the (open) left complex half plane, some require that the pencil $\lambda E-A$ is furthermore regular and of index at most one, since otherwise arbitrary small perturbations make the system unstable, see \cite{DuLM13,LinM09,MehU23} for detailed discussions, which also include the robustness question when the pencil $\lambda E-A$ is close to singular or high index.

In this paper we address the problem of determining {\color{black} proportional and/or derivative} output feedback controls 
that make the closed loop system regular and of index at most one, i.e. uniquely solvable for consistent initial conditions, and also  asymptotically stable. We study this problem for the  important class of port-Hamiltonian descriptor system representations that are introduced in the next subsection.

\subsection{Port-Hamiltonian descriptor systems}\label{sec:phDAE}
In this subsection we introduce the framework of port-Hamiltonian descriptor systems.
\begin{definition}\label{def:pHDAE}
A linear time-invariant descriptor system of the form
\begin{eqnarray}
E \dot x &=&(J-R) Qx + (B-P)u, \nonumber\\
y&=& (B+P)^H Q x  + (S-N) u, \label{genphdae}
\end{eqnarray}
with $E,Q\in {\mathbb C}^{\ell,n}$, $J,R\in {\mathbb C}^{\ell,\ell}$,
$B,P\in{\mathbb C}^{\ell,m}$,  $S=S^H, N=-N^H\in {\mathbb C}^{m,m}$
is called  \emph{port-Hamiltonian differential-algebraic (pHDAE) system with quadratic nonnegative Hamiltonian}
\begin{equation}\label{defham}
\mathcal H(x): = \frac 12 \Re(x^HQ^HEx) \geq 0
\end{equation}
if the following properties are satisfied:
\begin{itemize}
\item [i)]    $0\leq Q^HE = E^HQ\in \mathbb C^{n,n}$ and
$0=\Re(Q^H(J-J^H)Q)$;
\item [ii)] the \emph{dissipation matrix}
\begin{equation} \label{Wdef}
W=\left[\begin{array}{lc}
Q^H R Q& Q^H P \\[1mm]
P^H Q & S
\end{array}\right] \in \mathbb{C}^{n+m,n+m}
\end{equation}
is positive semidefinite, i.e., $W=W^H\geq 0$.
\end{itemize}
\end{definition}

The class of pHDAE  systems  provides a unified and natural modeling framework for the simulation and control of almost all classes of real world physical systems, see \cite{BeaMXZ18,JacZ12,MehM19,MehS23,MehU23,SchJ14,SchM23} for detailed discussions and a multitude of applications. The great success of modeling with pHDAE systems is mainly due to its many important properties.

Key properties of pHDAEs, see e.g. \cite{MehU23}, are the invariance of the class under power-conserving interconnection, which allows modularized automated modeling, the invariance under Galerkin projection which makes them ideal for discretization and model reduction, and in particular the encoding of properties like energy dissipation, stability and passivity in the algebraic structure of the coefficients of the equations.  The class of pHDAE systems also provides an ideal framework for  robust and physically interpretable control design. This follows, in particular, from the power balance equation and the resulting dissipation inequality,
see e.g. \cite{MehM19}.
\begin{theorem}\label{thm:pbe}
Consider a pHDAE system of the form \eqref{genphdae}. {\color{black}Then for any input $u$ the power balance equation 
  \begin{equation}\label{eq:powerBalanceEq}
    \frac{d}{dt}\mathcal H(x) = - \begin{bmatrix}
        x\\ u
    \end{bmatrix}^HW \begin{bmatrix}
        x\\ u
    \end{bmatrix}+ \Re (y^Hu)
  \end{equation}
  }
  holds along any solution $x$.
  In particular, the \emph{dissipation inequality}
  \begin{equation}\label{eq:dissIneq}
    \mathcal H(x(t_2)) - \mathcal H(x(t_1)) \leq \int_{t_1}^{t_2}\Re (y(\tau)^Hu(\tau))\ d\tau
  \end{equation}
  holds.
\end{theorem}

In physical space, one can view pHDAE systems as modeling the interaction of three types of energies by encoding these  in the structure of the coefficients. The \emph{stored energy} is presented by the nonnegative Hamiltonian $\mathcal H(x)$, the \emph{dissipated energy}  by the nonnegative quadratic form
$ \mathcal D(x,u)=
\begin{bmatrix}
    x \\ u
\end{bmatrix}^H W
\begin{bmatrix}
    x \\ u
\end{bmatrix}
$ and the \emph{supplied energy} by
$\mathcal S(y,u)=\Re( y^H u)$.

While for general descriptor systems it is computationally difficult to analyze whether a system is asymptotically stable, see e.g. \cite{BenBMX02,ZhoDG96}, in the pHDAE modeling framework, using Theorem~\ref{thm:pbe} easily allows to analyze when a pHDAE system is  stable (asymptotically stable).
It is well known  \cite{MehMW18} that if $Q$ has full column rank, then the pHDAE systems of the form \eqref{genphdae} are stable (but not necessarily asymptotically stable) in the sense that all finite eigenvalues are in the closed left complex half plane and those on the imaginary axis are semisimple. Furthermore, it is shown in  \cite{MehMW18} that the index of a pHDAE system can be at most $\nu=2$ and in \cite{MehMW21} the singularity is characterized by a common nullspace property.
Furthermore,
if the system is in the pHDAE representation, it is only needed to check the semidefiniteness of $E^HQ$ and $W$, which can be done accurately and with perturbation bounds via the calculation of Cholesky decompositions, see e.g. \cite{GolV96}.

There also exist structure preserving versions of the Kronecker canonical form, see \cite{AchAM21,BeaGM22}, where in order to preserve the  structure  and in particular the different types of energy $\mathcal H,\mathcal D,\mathcal S$, we require the transformations to satisfy $S=T^H$ and $Y=V^{-H}$, see \cite{BeaMX22,MehS23}.
We will discuss such condensed forms in Section~\ref{sec:condensed}.

{\color{black} It has been addressed in \cite{MehU23} how one can  reformulate a general linear pHDAE system to one with $\ell=n$ and $Q=I$ and how to remove the feedthrough term, so that $Du=(S-N)u=0$.
Although it will always be the first step of the regularization and stabilization procedure, we do not present this simplification here, but assume} that we have given a pHDAE system of the form
\begin{eqnarray}
E \dot x &=&(J-R)x + Bu, \nonumber\\
y&=& B^H  x, \label{phdae}
\end{eqnarray}
with $E,J,R\in {\mathbb C}^{n,n}$,
$B\in{\mathbb C}^{n,m}$,  $E=E^H\geq 0$, $R=R^H\geq 0$, $J=-J^H$, with the quadratic Hamiltonian
$\mathcal H(x)=\frac 12 x^H E x\geq 0$ and the dissipation matrix $W=\mat{cc} R & 0 \\ 0 & 0 \rix\geq 0$.
We also assume, without loss of generality, that $B$ has full column rank by restricting, if necessary, $u,y$ to an appropriate subspace.

\subsection{Problem statements}\label{S2}
{\color{black} For general unstructured descriptor systems the modification of system properties like regularity or stability via feedback  has been studied extensively, see e.g. \cite{BunBMN99,CamKM12,ChuCH98,ChuH99,NicC15,OzcL90,ShaZ87}.}
But such general feedback  approaches  do not necessarily preserve the pHDAE structure.
For pHDAE systems the  natural feedback classes  are proportional output feedbacks, since then the symmetry structure of the coefficients is preserved and it is sufficient if the feedback preserves the nonnegativity of the energy functions $\mathcal H$ and $\mathcal D$. We therefore discuss  proportional output feedback of the form
\[
u(t)=(F_S-F_H) y(t) +v(t)
\]
where $F_S=-F_S^H$ and $F_H=F_H^H$ are   such that  the resulting closed loop system
\begin{eqnarray*}
   \dot x(t) &=& (J+BF_SB^H -(R+BF_HB^H))x(t)+Bv(t),
   \nonumber\\
              y(t) &=& B^Hx(t),
\end{eqnarray*}
has desired properties.
In particular, we study the following three problems:

\textbf{Problem 1} (Regularization of pHDAE system (\ref{phdae})   by proportional output feedback):
Determine matrices  $F_S=-F_S^H, F_H=F_H^H$ such that
 the  pair $(E, J+BF_SB^H-(R+ BF_HB^H))$ is regular,
and $R+BF_HB^H\geq 0$, i.e.,  the resulting closed-loop system
%
%
is a regular pHDAE system.

\textbf{Problem 2} (Regularization and index reduction of pHDAE system (\ref{phdae})   by proportional output feedback):
Determine matrices  $F_S=-F_S^H, F_H=F_H^H$ such that
 the  pair $(E, J+BF_SB^H-(R+ BF_HB^H))$ is regular,  of index at most one,  and
$R+BF_HB^H\geq 0$, i.e.,  the resulting closed-loop system
is a regular pHDAE system of index at most one.

\textbf{Problem 3} (Stabilization of pHDAE system (\ref{phdae}) by proportional output feedback):
Determine  matrices  $F_S=-F_S ^H, F_H=F_H ^H$ such that
 the   pair $(E, J+BF_SB ^H-(R+ BF_HB ^H))$ is regular,  of index at most one, has all its finite eigenvalues in the open left {\color{black}  complex half plane}, and
 \[ R+BF_HB^H\geq 0, \]
 i.e.,  the resulting closed-loop system is a regular pHDAE system of index at most one and has all its finite eigenvalues in the open left {\color{black} complex half plane}.

In some applications it is also possible to use derivative output feedback $u=K\dot y$ to perform regularization, index reduction and stabilization. {\color{black} Our results   also extend to this case, see  Section~\ref{sec:derivative feedback}.}

All the constructions and conditions that we present are derived via structured condensed forms that we present in the next section. For completeness we also present coordinate free versions of the results for which we denote a full column rank matrix with its columns spanning the right  nullspace of a matrix $M$ by $\mathcal S_\infty(M)$ and with its columns spanning the left nullspace of $M$ by  $\mathcal T_\infty(M)$, respectively.
\section{Condensed forms}\label{sec:condensed}
The basis for the construction of regularizing feedbacks is the computation of condensed forms. In order to be able to construct the regularizing feedbacks in a numerically stable way we use unitary transformations. The following form is a modification of the condensed form presented in \cite{BeaGM22}.
\begin{lemma}\label{lemma1} Consider a pHDAE system of the form (\ref{phdae}).  Then there exist  unitary matrices $U$ and $V$ such that
\begin{small}
\begin{eqnarray}
U^H B V  &=& \bmat{ & m-n_3   & n_3 \cr
    n_1  & 0           & B_{12}    \cr
    n_2  &  B_{21} & B_{22}  \cr
    n_3  &  0         & B_{32}  \cr
    n_4  & 0          & 0           \cr
    n_5  & 0         & 0          \cr
    n_6  & 0         & 0         \cr},  \
U^H E U = \bmat{  & n_1 & n_2 & n_3 & n_4 & n_5 & n_6 \cr
    n_1 & E_{11}     & E_{12}      & E_{13}    & E_{14} & 0 & 0 \cr
    n_2 & E_{12}^H & E_{22}     & E_{23}     & E_{24} & 0 & 0 \cr
    n_3 & E_{13}^H & E_{23}^H & E_{33}     & E_{34} & 0 & 0 \cr
    n_4 & E_{14}^H & E_{24}^H & E_{34}^H & E_{44} & 0 & 0 \cr
    n_5 & 0              & 0          & 0        & 0  & 0  & 0 \cr
    n_6 & 0              & 0          & 0        & 0  & 0  & 0  \cr},   \label{condensed-form1}  \\
U^H(J-R)U&=&
\bmat{ & n_1 & n_2 & n_3 & n_4 & n_5  & n_6  \cr
    n_1  & J_{11}-R_{11}           & J_{12}-R_{12}           & J_{13}-R_{13}           & J_{14}-R_{14}        & J_{15}-R_{15}                    & J_{16} \cr
    n_2  &- J_{12}^H-R_{12}^H & J_{22}-R_{22}           & J_{23}-R_{23}           & J_{24}-R_{24}         & J_{25}-R_{25}                    & J_{26} \cr
    n_3  & -J_{13}^H-R_{13}^H & -J_{23}^H-R_{23}^H & J_{33}-R_{33}           & J_{34}-R_{34}          &  J_{35}-R_{35}                   & 0        \cr
   n_4  & -J_{14}^H-R_{14}^H & -J_{24}^H-R_{24}^H & -J_{34}^H-R_{34}^H & J_{44}-R_{44}          & J_{45}-R_{45}                    & 0        \cr
   n_5  & -J_{15}^H-R_{15}^H & -J_{25}^H-R_{25}^H & -J_{35}^H-R_{35}^H & -J_{45}^H-R_{45}^H & J_{55}-R_{55}                   & 0  \cr
   n_6  & -J_{16}^H                 & -J_{26}^H          & 0                               & 0                    & 0                    & 0 \cr}, \nonumber
\end{eqnarray}
\end{small}
where
\begin{equation}\label{eq1}  \rank \mat{cc} J_{16} \\ J_{26}\rix=n_1+n_2, \ \rank(B_{21})=n_2, \ \rank(B_{32})=n_3,  \ \rank(J_{55}-R_{55})=n_5,
\end{equation}
and, furthermore,
\begin{equation}\label{eq2} \rank \mat{cccccc} E_{11}     & E_{12}     & E_{13}      & E_{14} & 0         & B_{12} \\
 E_{12}^H & E_{22}     & E_{23}      & E_{24} & B_{21} & B_{22} \\
 E_{13}^H & E_{23}^H & E_{33}      & E_{34} & 0         & B_{32} \\
  E_{14}^H & E_{24}^H & E_{34}^H & E_{44} & 0         & 0 \rix=n_1+n_2+n_3+n_4, \  \ E_{44}>0.
\end{equation}
\end{lemma}
\proof A constructive proof that can be directly implemented as a numerical method is presented in  Appendix A.
\eproof

If  one allows nonunitary transformations in Lemma~\ref{lemma1}, then one can reduce the condensed form further.
\begin{corollary}\label{lemma2} Consider a pHDAE system of the form (\ref{phdae}).  Then there exist  nonsingular matrices $S$, $T$, and a unitary matrix $V$
such that
\begin{eqnarray}
S B V &=&\bmat{ & m-n_3   & n_3 \cr
n_1  & 0           & 0    \cr
n_2  &  B_{21} & 0  \cr
n_3  &  0         & B_{32}  \cr
n_4  & 0          & 0           \cr
n_5  & 0         & 0          \cr
n_6  & 0         & 0         \cr},  \
\ S E T = \bmat{  & n_1 & n_2 & n_3 & n_4 & n_5 & n_6 \cr
    n_1 & E_{11}     & 0                & E_{13}    & 0 & 0 & 0 \cr
    n_2 & 0                & E_{22}      & E_{23}     & 0 & 0 & 0 \cr
    n_3 & E_{13}^H & E_{23}^H & E_{33}     & 0 & 0 & 0 \cr
    n_4 & 0 & 0 & 0 & E_{44} & 0 & 0 \cr
    n_5 & 0              & 0          & 0        & 0  & 0  & 0 \cr
    n_6 & 0              & 0          & 0        & 0  & 0  & 0  \cr} 
    \nonumber \\
S(J-R)T &=&  \bmat{ & n_1 & n_2 & n_3 & n_4 & n_5  & n_6  \cr
    n_1  & A_{11}& A_{12}   & A_{13}   & A_{14}        & 0                    & A_{16} \cr
    n_2  & A_{21} & A_{22}  & A_{23}   & A_{24}        & 0                    & A_{26} \cr
    n_3  & A_{31} & A_{32}  & A_{33}   & A_{34}        & 0                    & 0        \cr
    n_4  & A_{41} & A_{42}  & A_{43}   & A_{44}        & 0                    & 0        \cr
    n_5  & 0 & 0 & 0 & 0 & A_{55}                   & 0  \cr
    n_6  & -A_{16}^H                 & -A_{26}^H          & 0                               & 0                    & 0                    & 0 \cr}, \label{condensed-form2}
\end{eqnarray}
where $SB=T^HB$,
\begin{equation}\label{eq3}  \rank \mat{cc} A_{16} \\ A_{26}\rix=n_1+n_2, \quad \rank(B_{21})=n_2, \quad \rank(B_{32})=n_3,  \quad \rank(A_{55})=n_5,
\end{equation}
and
\begin{equation}\label{eq4} E_{11}> 0, \ E_{44}> 0,
\mat{ccc}              E_{11}     & 0                & E_{13}   \\
   0             & E_{22}     & E_{23}   \\
   E_{13}^H & E_{23}^H & E_{33}  \rix \geq 0.
\end{equation}
\end{corollary}
\proof The proof follows by block Gaussian elimination in \eqref{condensed-form1}.
\eproof
Using Corollary~\ref{lemma2} we immediately obtain the following coordinate-free descriptions of the dimensions in the condensed form~\eqref{condensed-form1}.
\begin{corollary}\label{lemma3} Consider a pHDAE system of the form \eqref{phdae} in the condensed form (\ref{condensed-form2}). Then the following statements hold.

i)
\begin{eqnarray*} n_1+n_4&=&\rank \mat{cc} E & B \rix-\rank(B),   \\
n_3+n_4&=&\rank( \cT_\infty^H((J-R)\cS_\infty(\mat{c} E \\ B^H \rix) ) \mat{cc} E & B \rix).
\end{eqnarray*}

ii)
\begin{equation} \rank \mat{ccc} E & J-R & B \rix=n \label{cond1}
\end{equation}
if and only if
\begin{equation} n_6=n_1+n_2.  \label{cond2}
\end{equation}

iii)
\[  \rank(E_{13})=\rank( \cT^H_\infty(B) E \cS_\infty( \cT^H_\infty(\mat{cc} E & B \rix)(J-R))-n_4, \]
and
$\rank(E_{13})=n_1 $
if and only if
\begin{equation} \rank(\cT^H_\infty(B) E \cS_\infty( \cT^H_\infty(\mat{cc} E & B \rix)(J-R)  ) )=\rank \mat{cc} E & B \rix-\rank(B).  \label{cond3}
\end{equation}
\end{corollary}

\proof   {\color{black} To read off the related spaces $\cT_\infty$ and $\cS_\infty$ from  the condensed form in Corollary~\ref{lemma2}, we assume without loss of the generality that
$S=I$, $T=I$ and $V=I$ in Corollary~\ref{lemma2}.

 i) We have
\[
\rank \mat{cc} E & B \rix-\rank(B)=(n_1+n_2+n_3+n_4)-(n_2+n_3)=n_1+n_4.
\]
It then follows that
 \begin{eqnarray*}
  \rank( \cT_\infty^H((J-R)\cS_\infty  \left (\mat{c} E \\ B^H \rix\right ))  \mat{cc} E & B \rix)
&=& \rank ( \cT_\infty^H(\mat{cc} 0 & A_{16} \\ 0 & A_{26} \\ 0 & 0 \\ 0 & 0 \\ A_{55} & 0 \\ 0 & 0 \rix)  \mat{cc} E & B \rix) \\
&=& \rank (
\mat{ccccc} E_{13}^H & E_{23}^H & E_{33} & 0      & B_{32} \\
 0             & 0              & 0     & E_{44} & 0      \rix)   \\
 &=&  n_3+n_4.
\end{eqnarray*}

ii) Since
\begin{eqnarray*}
  \rank \mat{ccc} E & J-R & B \rix &=& n_1+n_2+n_3+n_4+n_5+\rank \mat{cc} -A_{16}^H & -A_{26}^H \rix \\
    &=& n-n_6+(n_1+n_2),
\end{eqnarray*}
it follows that  (\ref{cond1}) holds if and only if (\ref{cond2}) holds.

 iii)  Note that
\[ \cT^H_\infty(B) E=\mat{cccccc} E_{11} & 0 & E_{13}  & 0 & 0 & 0 \\          0 & 0 & 0 & E_{44} & 0 & 0 \\
0 & 0 & 0 & 0 & 0 & 0 \\
0 & 0 & 0 & 0 & 0 & 0 \rix,
\]
and
\[ \cT^H_\infty(\mat{cc} E & B \rix)(J-R)
    =\mat{cccccc} 0 & 0 & 0 & 0 & A_{55} & 0 \\
-A_{16}^H & -A_{26}^T & 0 & 0 & 0 & 0 \rix. \]
Hence, we obtain
\begin{eqnarray*}
\rank( \cT^H_\infty(B) E \cS_\infty( \cT^H_\infty(\mat{cc} E & B \rix)(J-R))-n_4
&=& \rank \mat{ccc} E_{13} & 0 & 0 \\
0 & E_{44} & 0 \\
0 & 0 & 0 \\
0 & 0 & 0 \rix-n_4 \\
&=& \rank(E_{13}).
\end{eqnarray*}
Furthermore, we have $\rank(E_{13})=n_1 $ if and only if
\begin{eqnarray*}
 \rank( \cT^H_\infty(B) E \cS_\infty( \cT^H_\infty(\mat{cc} E & B \rix)(J-R)) &=& n_1+n_4=(n_1+n_2+n_3+n_4)-(n_2+n_3) \\
  &=& \rank \mat{cc} E & B \rix-\rank(B).
\end{eqnarray*}
}
\hskip 12cm
\eproof

The condensed forms in this section form the basis for the solution of  problems 1-3 in the following section.

\section{Regularization and stabilization via proportional output feedback}\label{sec:regulfeed}
In this section we characterize the solutions of  Problems 1--3.
The characterizations of the solution to the first two problems have similar conditions as  in the unstructured case.
\begin{theorem}\label{theorem1} Consider a pHDAE system (\ref{phdae}).  Then Problem 1 
 is solvable  if and only if  (\ref{cond1}) holds.
\end{theorem}
\proof   Suppose that there exist matrices  $F_S=-F_S^H$ and $F_H=F_H^H$  such that $(E, J+BF_SB^H-(R+BF_HB^H)$ is regular. Then we have
\[
\det (sE-(J+BF_SB^H-(R+BF_HB^H)))\not=0, \quad {\rm for \ some \ } s\in \mathbb C,
\]
which together with the condensed form (\ref{condensed-form2}) gives the condition (\ref{cond2}). Then, Corollary~\ref{lemma3} ii) yields  the condition (\ref{cond1}). Hence, the necessity follows.

To show the sufficiency, let the condition (\ref{cond1}) and thus equivalently (\ref{cond2}) holds.  Let $F_{22}\in \mathbb C^{n_3, n_3}$ be such that
$F_{22}>0$ and that $A_{33}-B_{32}F_{22}B_{32}^H$ is nonsingular. {\color{black} This is possible since $B_{32}$ has full row rank}. Then with
\[F_S=0, \quad F_H=V \mat{cc} 0 & 0 \\ 0 & F_{22}  \rix V^H, \]
we have that
$(E, J+BF_SB^H-(R+BF_HB^H))$ is regular and $R+BF_HB^H\geq 0$.   \eproof

If we further require that the index of the closed loop system pencil is reduced to one then we have the following result.
\begin{theorem}\label{theorem4} Consider a pHDAE system of the form (\ref{phdae}).  Then Problem 2  is solvable
 if and only if
\begin{equation}\label{con-1}
               \rank \mat{ccc} E & (J-R) \cS_\infty(E) & B \rix=n.
\end{equation}
\end{theorem}

\proof
 {\color{black}  Let  $F=F_S-F_H$, with $F_S=-F_S^H$ and $F_H=F_H^H$, be such that $(E, J-R+BFB^H)$ is regular and of index at most one.
Set
\[ V^HFV=\bmat{ & m-n_3 & n_3 \cr
             m-n_3     & F_{11} & F_{12} \cr
                 n_3     & F_{21} & F_{22} \cr}. \]
Then condition (\ref{cond1}) holds  and,
denoting by ${\rm deg}(p(s)$ the degree of the polynomial $p(s)$, we have
\[ {\rm deg}\det (sE-(J-R+BFB^T))=\rank(E),  \]
i.e.,
\begin{equation}\label{R1}  n_1+n_2=n_6, \quad \mat{c} A_{16} \\ A_{26} \rix {\rm \  is \ nonsingular},
\end{equation}
 and
\begin{eqnarray*}
 {\rm deg}\det (sE-(J-R+BFB^T)) &=&  {\rm deg}\det (S(sE-(J-R+BFB^T)T) \\
&=& {\rm deg}\det \left (\mat{cc} sE_{33}-(A_{33}+B_{32}F_{22}B_{32}^H) & -A_{34} \\ -A_{43} & sE_{44}-A_{44} \rix \right ) \\
&=& \rank(E_{33})+\rank(E_{44}) \\
&=&\rank(E) \\
&=& \rank \mat{ccc} E_{11} & 0 & E_{13} \\ 0 & E_{22} & E_{23} \\ E_{13}^T & E_{23}^H & E_{33} \rix+\rank(E_{44}),
\end{eqnarray*}
which gives
\begin{equation}\label{cond6} \rank \mat{ccc} E_{11} & 0 & E_{13} \\ 0 & E_{22} & E_{23} \\ E_{13}^H & E_{23}^H & E_{33} \rix
              =\rank (E_{33}).
\end{equation}
Note that  $\mat{ccc}              E_{11}     & 0                & E_{13}   \\
   0             & E_{22}     & E_{23}   \\
   E_{13}^H & E_{23}^H & E_{33}  \rix \geq 0$,
and thus, condition (\ref{cond6}) is equivalent to
\begin{equation}\label{cond7}
   \mat{cc} E_{11} & 0 \\ 0 & E_{22} \rix -\mat{c} E_{13} \\ E_{23} \rix E_{33}^{+}\mat{cc} E_{13}^H & E_{23}^H \rix=0,
\end{equation}
where $E_{33}^{+}$ is the Moore-Penrose inverse of $E_{33}$.
A direct calculation yields that conditions (\ref{R1}) and (\ref{cond7}) imply condition (\ref{con-1}). Hence, the necessity follows.   }

To show the sufficiency, take
\[
F_S=0, \quad F_H=V\mat{cc} 0 & 0  \\ 0  & F_{22} \rix V^H, \]
with  $F_{22}>0$,  and $  \cT_\infty^H(E_{33})(A_{33}-B_{32}F_{22}B_{32}^H)\cS_\infty(E_{33})$ nonsingular.
Then the pair $(E_{33}, A_{33}-B_{32}F_{22}B_{32}^H)$  is regular and of index at most one. Because  the condition (\ref{con-1}) implies the conditions (\ref{cond1})
and (\ref{cond7}), we have that $ R+BF_HB^H \geq 0$ and the pair $(E, J+BF_SB^H-(R+BF_HB^H))$ is regular and of index at most one.
\eproof

After  a pHDAE system  of the form (\ref{phdae}) has been regularized and made of index at most one, the next task is to design a proportional output feedback so that the resulting closed-loop system is
\emph{asymptotically stable}, i.e.  all its finite eigenvalues have negative real part. While this a very hard and partially unsolved problem for general descriptor systems, for pHDAE systems the solution is surprisingly simple.

We need the following lemma.
\begin{lemma} \label{lemma-S1} Consider $E, J, R\in \mathbb C^{n,n}$ with $E\geq 0$, $J=-J ^H$, $R\geq 0$, and
\[
E=\bmat{ & n_1 & n_2 \cr
             n_1 & E_{11} & E_{12} \cr
             n_2 & E_{12} ^H & E_{22} \cr}, \quad
   R=\bmat{ & n_1 & n_2 \cr
             n_1& R_{11} & 0 \cr
             n_2 & 0         & 0   \cr}, \quad
   J=\bmat{ & n_1 & {\color{black} n_2} \cr
             n_1 & J_{11} & J_{12} \cr
             n_2 & -J_{12} ^H & J_{22} \cr}, \]
where $R_{11}>0$. Then the following statements hold.
\begin{itemize}
\item[ i)]  $J-R$ is nonsingular if and only if
\[ \rank \mat{c} J_{12} \\ J_{22} \rix=n_2. \]
\item[ii) ]  The  pair   $(E, J-R)$ has all its finite eigenvalues in the open left complex half plane if and only if for all purely imaginary $s\in \mathbb C$
\begin{equation}\label{stable-1}
\rank \mat{c} J_{12}-sE_{12} \\
    J_{22}-sE_{22}  \rix=n_2.
\end{equation}
\end{itemize}
\end{lemma}
\proof The proof of i) is trivial.

ii) If the  pair   $(E, J-R)$ has all its finite eigenvalues in the open left complex half plane, then obviously (\ref{stable-1}) holds for all purely imaginary $s\in \mathbb C$.

Conversely,  let (\ref{stable-1}) hold for all purely imaginary $s\in \mathbb C$.
It follows from i) that $(E, J-R)$ is regular.  Next, let $s_0\in \mathbb C$ be any finite eigenvalue of $(E, J-R)$, and let $x=\begin{bmatrix} x_1\\
 x_2 \end{bmatrix} \in \mathbb C^n$ (partitioned analogously) be a corresponding eigenvector
normalized such that $x^H Ex=1$. Then we have
\begin{equation} \label{eig-1}
       \mat{cc} J_{11}-R_{11} & J_{12} \\ -J_{12} ^H & J_{22} \rix x=s_0 \mat{cc} E_{11} & E_{12} \\ E_{12} ^H & E_{22} \rix x,
\end{equation}
and hence,
\begin{eqnarray*}
  s_0 &=& {\color{black} s_0 } ~x ^H  \mat{cc} E_{11} & E_{12} \\ E_{12} ^H & E_{22} \rix x \\
         &=& x ^H  \mat{cc} J_{11}-R_{11} & J_{12} \\ -J_{12} ^H & J_{22} \rix x
                  =-x_1 ^H R_{11}x_1+x ^H  \mat{cc} J_{11} & J_{12} \\ -J_{12} ^H & J_{22} \rix x,
\end{eqnarray*}
which gives
\[
\Re(s_0)=-x_1 ^H R_{11}x_1\leq 0.
\]
We show that $x_1\not=0$. If we had $x_1=0$, then $s_0$ is purely imaginary, $x_2\not=0$ and
\[ \mat{c} J_{12}-s_0E_{12} \\
                       J_{22}-s_0E_{22}  \rix x_2=0.
\]
This and the condition that $\rank \mat{c} J_{12}-s_0E_{12} \\
    J_{22}-s_0E_{22}  \rix=n_2$ yields that $x_2=0$ which is a contradiction. Hence, $x_1\not=0$ and
$\Re(s_0)=-x_1^H R_{11}x_1<0$. Therefore, $(E, J-R)$ has all its finite eigenvalues in the open left complex half plane.
\eproof

%

We now present necessary and sufficient solvability conditions for the solution of Problem 3.
\begin{theorem}\label{stability-1}  Consider  a pHDAE system of the form (\ref{phdae}).  Then Problem 3 
 is solvable  if and only if the condition (\ref{con-1}) holds and for all purely imaginary $s\in \mathbb C$
\begin{equation}\label{con-S1}
\rank \mat{cc} J-R-sE & B \rix=n.
\end{equation}
\end{theorem}

\proof  {\color{black} Let the matrix $F=F_S-F_H$, with $F_S=-F_S^H$ and $F_H=F_H^H$, be such that  the   pair $(E, J-R+ BFB^H )$ is regular,  of index at most one, has all its finite eigenvalues in the open left {\color{black}  complex half plane}, and
 $R+BF_HB^H\geq 0$. Then  condition (\ref{con-1}) follows by Theorem~\ref{theorem4}, and moreover, for all purely imaginary $s\in \mathbb C$,
\[ \rank \mat{cc} J-R-sE & B \rix=\rank \mat{cc} J-R+BFB^H-sE & B \rix=n. \]
Hence, the necessity  follows.  }

To prove the sufficiency, let $U$ be a unitary matrix such that
\[
U^HB=\bmat{    &    \cr
    n_1 & 0 \cr
    n_2 &  B_2 \cr
    n_3 &  0 \cr}, \quad
 U^HRU=\bmat{  & n_1 & n_2 & n_3 \cr
 n_1   & R_{11} & R_{12} & 0 \cr
 n_2   &  R_{12} ^H & R_{22} & 0 \cr
 n_3   &  0 & 0 & 0 \cr},
\]
where
\[
\rank(B_2)=n_2, \ R_{11}>0.
\]
Set
\[
U^HJU=\bmat{ & n_1 & n_2 & n_3 \cr
     n_1 & J_{11} & J_{12} & J_{13} \cr
    n_2 & -J_{12} ^H & J_{22} & J_{23} \cr
    n_3 &  -J_{13} ^H & -J_{23} ^H & J_{33} \cr}, \quad
   U^HEU=\bmat{ & n_1 & n_2 & n_3 \cr
    n_1 &  E_{11} & E_{12} & E_{13} \cr
    n_2 & E_{12} ^H & E_{22} & E_{23} \cr
    n_3 &  E_{13} ^H & E_{23} ^H & E_{33} \cr}.
\]
Let $F_S=0$ and $F_H=F_H^H$ be such that
\[
\rank(R+BF_HB ^H)=\rank\mat{cc} R & B \rix, \quad R+BFB ^H\geq 0,
\]
i.e.,
\[  \mat{cc} R_{11} & R_{12} \\ R_{12} ^H & R_{22}+B_2F_HB_2 ^H\rix >0.
\]
Then it follows from Lemma \ref{lemma-S1} and the fact that (\ref{con-S1}) holds for all purely imaginary $s$ that the  pair $(E, J-(R+BF_HB ^H))$ has all its finite eigenvalues in the open left complex half plane.

Next, let  $\tilde U\in \mathbb C^{n,n}$ be unitary such that
\[
\tilde U^H E \tilde U=\bmat{ & \tau_1 & \tau_2 & \tau_3 \cr
 \tau_1     & \tilde E_{11} & 0 & 0 \cr
\tau_2      & 0           & 0 & 0 \cr
 \tau_3      & 0            & 0 & 0 \cr}, \quad
    \tilde U^H (R+BF_HB^H) \tilde U
    =\bmat{ & \tau_1 & \tau_2 & \tau_3 \cr
 \tau_1   & \tilde R_{11} & \tilde R_{12} & 0   \cr
\tau_2   & \tilde R_{12}^H & \tilde R_{22} & 0 \cr
\tau_3    & 0              & 0          & 0 \cr},
\]
where
\[ \tilde E_{11}>0, \quad \tilde R_{22}>0.
\]
Set
\[ \tilde U^H J \tilde U=\bmat{ & \tau_1 & \tau_2 & \tau_3 \cr
\tau_1   & \tilde J_{11} & \tilde J_{12} & \tilde J_{13}   \cr
   \tau_2   & -\tilde J_{12}^H & \tilde J_{22} & \tilde J_{23}\cr
  \tau_3    & -\tilde J_{13}^H & -\tilde J_{23}^H  & \tilde J_{33} \cr}, \quad
 \tilde U^H B=\bmat{ &         \cr
  \tau_1 & \tilde B_1 \cr
\tau_2 & \tilde B_2 \cr
     \tau_3 & \tilde B_3 \cr}.
\]
Note that
\[
\rank(R+BF_HB ^H)=\rank\mat{cc} R & B \rix=\rank \mat{cc} R+BF_HB^H & B \rix, \]
and thus
\[ \tilde B_3=0.
\]
Additionally, condition  (\ref{con-1}) implies that
\[
\rank \mat{c} \tilde J_{23} \\ \tilde J_{33} \rix=\rank \mat{cc} -\tilde J_{23} ^H & \tilde J_{33} \rix=\tau_3,
\]
and hence by Lemma \ref{lemma-S1} we have that $\mat{cc} \tilde J_{22}-\tilde R_{22} & \tilde J_{23} \\ -\tilde J_{23}^H & \tilde J_{33} \rix$ is nonsingular. Therefore,
the  pair $(E, J-(R+BF_HB ^H))$ is regular and of index at most one.
\eproof

After the characterization of the existence of output feedbacks that make the pHDAE system regular and of index at most one as well as asymptotically stable an important question is to use the feedbacks in such a way that the resulting closed loop system is robustly regular, of index at most one and asymptotically stable.
In order to do this one needs efficiently computable characterizations what the distance to the nearest non-regular pHDAE, higher index pHDAE are
\cite{GilMS18,GugM22,MehMW21}, respectively the distance to instability \cite{AliMM20,GilS17,GilMS18} are. Furthermore, it is necessary to analyze how the pHDAE structure can be exploited, and how to compute robust pHDAE representations, see
\cite{BanMNV20,MehV20}.

\section{Regularization and stabilization via derivative output feedback}\label{sec:derivative feedback}
{\color{black} The results in the previous section can be generalized to the case that one includes also derivative feedback. Since derivative feedback is rarely used in practice, the following results are  interesting mainly from a theoretical point of view.}

\begin{theorem}\label{theorem2} Consider a pHDAE system  of the form (\ref{phdae}). There exists a derivative feedback matrix $K$  such that
 the pair $(E+BKB^H, J-R)$ is regular and $E+BKB^H\geq 0$  if and only if  (\ref{cond1}) holds.
\end{theorem}
\proof  Suppose there exists matrix $K$ such that $(E+BKB^H,  J-R)$ is regular. Then it follows that
\[ \det (s(E+BKB^H)-(J-R))\not=0, \quad {\rm for \ some \ } s\in \mathbb C,
\]
which together with  the condensed form  (\ref{condensed-form2}) gives  condition (\ref{cond2}), i.e. by equivalence also condition (\ref{cond1}) holds. Hence, the necessity is shown.

To show sufficiency, let $K_{22}\in \mathbb C^{n_3, n_3}$ be such that
$K_{22}>0$ and $E_{33}+B_{32}K_{22}B_{32}^H>0$. Taking
\[
K=V \mat{cc} 0 & 0 \\ 0 & K_{22}  \rix V^H,
\]
it follows  from (\ref{cond1}) (or equivalently from (\ref{cond2}))  that $(E+BKB^H, J-R)$ is regular and $E+BKB^H\geq 0$. Hence, the sufficiency is proved.
\eproof

We can also combine Theorems~\ref{theorem1} and~\ref{theorem2}.
\begin{theorem}\label{theorem3} Consider a pHDAE system of the form (\ref{phdae}).  There exist feedback matrices $K$, $F_S=-F_S^H$ and $F_H=F_H^H$ such that   the  pair $(E+BKB^H, J+BF_SB^H-(R+BF_HB^H)$ is regular and
$E+BKB^H\geq 0$, $R+BF_HB^H\geq 0$ if and only if  (\ref{cond1}) holds.  Moreover, if the condition (\ref{cond1}) holds, then for any integer $r$ satisfying
\begin{equation} \rank \mat{cc} E & B \rix-\rank(B) \leq r \leq \rank \mat{cc} E & B \rix, \label{cond4}
\end{equation}
there exist matrices $K$ and $F_H=F_H^H$ and $F_S=0$  such that $(E+BKB^H, J+BF_SB^H-(R+BF_HB^H))$ is regular, and
\[  \rank(E+BKB^H)=r, \quad E+BKB^H\geq 0, \quad R+BF_HB^H\geq 0.   \]
\end{theorem}

\proof  Suppose that there exist matrices $K$,  $F_S=-F_S^H$, and $F_H=F_H^H$  such that $(E+BKB^H, J+BF_SB^H-(R+BF_HB^H)$ is regular. Then
\[
\det (sE-(J+BF_SB^H-(R+BF_HB^H)))\not=0, \quad {\rm for \ some \ } s\in \mathbb C,
\]
from which we obtain condition (\ref{cond2}), and equivalently (\ref{cond1}). Hence, the necessity is shown.

By Corollary~\ref{lemma3},  conditions (\ref{cond1}) and (\ref{cond4}) are equivalent to condition (\ref{cond2}) and
\[
n_1+n_4\leq r \leq n_1+n_2+n_3+n_4,
\]
respectively.  Since $E_{11}>0$, $\rank(B_{21})=n_2$, $\rank(B_{32})=n_2$ and $E\geq 0$, there exists a matrix $K=\mat{cc} K_{11} & K_{12} \\ K_{12}^H & K_{22} \rix$ such that
\[
\mat{ccc} E_{11} & 0 & E_{13} \\ 0 & E_{22} & E_{23} \\ E_{13}^H & E_{23}^H & E_{33} \rix
   +\mat{cc} 0 & 0 \\ B_{21} & 0 \\ 0 & B_{32} \rix K \mat{cc} 0 & 0 \\ B_{21} & 0 \\ 0 & B_{32} \rix^H \geq 0, \]
\[  \rank\left (\mat{ccc} E_{11} & 0 & E_{13} \\ 0 & E_{22} & E_{23} \\ E_{13}^H & E_{23}^H & E_{33} \rix
   +\mat{cc} 0 & 0 \\ B_{21} & 0 \\ 0 & B_{32} \rix K \mat{cc} 0 & 0 \\ B_{21} & 0 \\ 0 & B_{32} \rix^H \right )=r-n_4. \]
Let $F_{22}>0$ be such that
\[
\cT_\infty^H(E_{33}+B_{32}K_{22}B_{32}^H) (A_{33}-B_{32}F_{22}B_{32}^H)
\cS_\infty(E_{33}+B_{32}K_{22}B_{32}^H)
\]
is nonsingular and set
\[
K=VK V^H,  \quad  F_H=V \mat{cc} 0 & 0 \\ 0 & F_{22}  \rix V^H, \quad F_S=0.
\]
We then have
\[
\rank(E+BKB^H)=r, \quad E+BKB^H\geq 0, \quad R+BF_HB^H\geq 0,
\]
and  $(E+BKB^H, J+BF_SB^H-(R+BF_HB^H))$ is regular.
\eproof

The corresponding results to achieve an index at most one are as follows.
\begin{theorem}\label{theorem5}  Consider a pHDAE system of the form (\ref{phdae}). There exists a matrix $K$  such that
 the  pair $(E+BKB^H, J-R)$ is regular and of index at most one, and $E+BKB^H\geq 0$ if and only if conditions (\ref{cond1}) and (\ref{cond3}) hold.
\end{theorem}
\proof  By Lemma \ref{lemma3},  conditions (\ref{cond1}) and (\ref{cond3}) are equivalent to
\[ n_6=n_1+n_2, \quad   \rank(E_{13})=n_1.
\]
If the pair $(E+BKB^H, J-R)$ is regular and of index at most one for some $K$, then with
\[  \mat{cc} K_{11} & K_{12} \\ K_{21} & K_{22} \rix =V^HKV,
\]
we have  condition (\ref{cond1}) and
\begin{equation}\label{cond8}
     \rank \left (\mat{ccc} E_{11} & 0 & E_{13} \\
   0           & E_{22}+B_{21}K_{11}B_{21}^H & E_{23}+B_{21}K_{12}B_{32}^H \\
E_{13}^H & E_{23}^H+B_{32}K_{21}B_{21}^H & E_{33}+B_{32}K_{22}B_{32}^H \rix \right )
                = \rank (E_{33}+B_{32}K_{22}B_{32}^H).
\end{equation}
It is obvious that (\ref{cond8}) implies
\[
E_{11}=E_{13}(  E_{33}+B_{32}K_{22}B_{32}^H )^{+} E_{13}^H,
\]
which together with $E_{11}>0$ gives $\rank(E_{13})=n_1$, i.e.,  condition (\ref{cond3}) holds. Hence, the necessity is shown.

To show the sufficiency, note that $E_{11}>0$, $\rank(E_{13})=n_1$, and $B_{32}$ is nonsingular, so there exists $K_{22}=K_{22}^H$ such that
\[
E_{33}+B_{32}K_{22} B_{32}^H> 0,
\]
and
\[ \rank \left( \mat{cc} E_{11} & E_{13} \\ E_{13}^H &  E_{33}+B_{32}K_{22} B_{32}^H \rix\right )=\rank( E_{33}+B_{32}K_{22} B_{32}^H)=n_3.
\]
Additionally, since $B_{21}$ is of full row rank, there exist $K_{11}$, $K_{12}$ such that
\[
E_{22}+B_{21} K_{11} B_{21}^H=0, \quad E_{23}+B_{21} K_{12} B_{32}^H=0.
\]
Taking
\[
K=V \mat{cc} K_{11} & K_{12} \\ K_{12}^H & K_{22} \rix V^H,
\]
we have that
\[ E+BKB^H\geq 0, \]
and $(E+BKB^H, J-R)$ is regular and of index at most one. \eproof

\begin{theorem}\label{theorem6}  Consider a pHDAE system of the form~(\ref{phdae}). There exist matrices $K$, $F_S=-F_S^H$, and $F_H=F_H^H$ such that  the  pair $(E+BKB^H, J+BF_SB^H-(R+BF_HB^H)$ is regular and of index at most one, and
$E+BKB^H\geq 0$, $R+BF_HB^H\geq 0$ if and only if conditions (\ref{cond1}) and (\ref{cond3}) hold. Moreover, under conditions (\ref{cond1}) and (\ref{cond3}),
for a given integer $r$, there exist matrices $K$,
$F_S=-F_S^H$ and $F_H=F_H^H$ such that
\[ E+BKB^H\geq 0, \quad R+BF_HB^H\geq 0, \]
$(E+BKB^H, (J+BF_SB^H)-(R+BF_HB^H))$ is regular,
$(E+BKB^H, (J+BF_SB^H)-(R+BF_HB^H))$ has index at most one
and $\rank(E+BKB^H)=r$
if and only if
\begin{equation}\label{cond9}
    \rank \mat{cc} E & B \rix-\rank(B)\leq r\leq \rank( \cT_\infty^H((J-R)\cS_\infty(\mat{c} E \\ B^H \rix) ) \mat{cc} E & B \rix).
\end{equation}
\end{theorem}

\proof  For any $K$ and $F$ with
\begin{equation}\label{KF} K=V \mat{cc} K_{11} & K_{12} \\ K_{21} & K_{22} \rix V^H, \quad
   F=V\mat{cc} F_{11} & F_{12} \\ F_{21} & F_{22} \rix V^H,
\end{equation}
it follows from direct calculation that  $(E+BKB^H, J-R+BFB^H)$ is regular and of index at most one
if and only if condition (\ref{cond1}) holds,
\begin{equation}\label{eq8} \rank \left (\mat{ccc} E_{11} & 0 & E_{13} \\
    0 & E_{22}+B_{21}K_{11}B_{21}^H & E_{23}+B_{21}K_{12}B_{32}^H \\
        E_{13}^H & E_{23}^H+B_{32}^H K_{21} B_{21}^H & E_{33}+B_{32}K_{22}B_{32}^H \rix\right )
   =\rank(E_{33}+B_{32}K_{22}B_{32}^H),
\end{equation}
 and $(E_{33}+B_{32}K_{22} B_{32}^H, A_{33}+B_{32}F_{22}\cB_{32}^H)$ is regular and of index at most one.
Obviously, (\ref{eq8}) with $E_{11}>0$ implies $\rank(E_{13})=n_1$, i.e., the condition (\ref{cond3}) holds.   Hence, necessity follows.

The sufficiency follows from the sufficiency of Theorem \ref{theorem5} with $F_S=0$ and $F_H=0$.

To study the possible rank of $E+BKB^H$, for any $K$ and $F$ of the form (\ref{KF}) with $(E+BKB^H, J-R+BFB^H)$ being regular and of index at most one, we obtain
\begin{eqnarray*}
 n_1+n_4 &\leq&   \rank \left( \mat{ccc} E_{11} & 0 & E_{13} \\
        0 & E_{22}+B_{21}K_{11}B_{21}^H & E_{23}+B_{21}K_{12}B_{32}^H \\
        E_{13}^H & E_{23}^H+B_{32}^H K_{21} B_{21}^H & E_{33}+B_{32}K_{22}B_{32}^H \rix\right ) +\rank(E_{44})  \\
               &=& \rank(E+BKB^H) \\
               &=& \rank (E_{33}+B_{32} K_{22} B_{32}^H)+\rank(E_{44}) \\
               &\leq& n_3+n_4,
\end{eqnarray*}
which together with Corollary \ref{lemma3} gives condition (\ref{cond9}).

Let $r$ be any integer satisfying the condition (\ref{cond9}). We can assume without loss of generality that
\[
E_{13}=\bmat{ &  n_1                & n_3-n_1 \cr
    & E_{13}^{(1)} & 0     \cr},
    \quad B_{32}=\bmat{ & n_1 & n_3-n_1 \cr
n_1  & B_{32}^{(1)}  & 0      \cr
n_3-n_1 & 0         & B_{32}^{(4)} \cr}, \]
where
\[
\rank( E_{13}^{(1)})=n_1, \quad \rank(B_{32}^{(1)})=n_1, \quad \rank(B_{32}^{(4)})=n_3-n_1.
\]
Set
\[ E_{33}=\bmat{ & n_1    &  n_3-n_1       \cr
    n_1  & E_{33}^{(1)}       & E_{33}^{(2)} \cr
    n_3-n_1  & (E_{33}^{(2)})^H & E_{33}^{(4)} \cr}, \quad
 A_{33}=\bmat{ & n_1    &  n_3-n_1       \cr
n_1  & A_{33}^{(1)}    & A_{33}^{(2)} \cr
     n_3-n_1  & A_{33}^{(3)} & A_{33}^{(4)} \cr}.
\]
Let $K_{11}$, $K_{12}$,  $K_{22}^{(1)}$, $K_{22}^{(2)}$ and $K_{22}^{(4)}$  be such that
\[ E_{22}+B_{21} K_{11} B_{21}^H=0, \quad E_{23}+B_{21} K_{12} B_{32}^H=0,
\]
\[ E_{33}^{(1)}+B_{32}^{(1)}K_{22}^{(1)}(B_{32}^{(1)})^H=(E_{13}^{(1)})^H E_{11}^{-1}E_{13}^{(1)}, \quad
 E_{32}^{(2)}+B_{32}^{(1)}K_{22}^{(2)}(B_{32}^{(4)})^H=0,
 \]
\[ E_{33}^{(4)}+B_{32}^{(4)}K_{22}^{(4)}(B_{32}^{(4)})^H=\mat{cc} \Lambda  & 0 \\ 0 & 0 \rix,
\]
and
\[
K_{22}=\mat{cc} K_{22}^{(1)} & K_{22}^{(2)} \\ (K_{22}^{(2)})^H & K_{22}^{(4)} \rix,
\]
where $\Lambda \in \mathbb C^{(r-n_1-n_4), (r-n_1-n_4)}$,  $\Lambda>0$. Furthermore, let
$F_{22}^{(4)} \in \mathbb C^{(n_1+n_3+n_4-r), (n_1+n_3+n_4-r)}$ satisfy that
\[
F_{22}^{(4)}>0, \quad  A_{33}^{(4)}-B_{32}^{(4)} F_{22}^{(4)} (B_{32}^{(4)})^H=\mat{cc} \star  & \star \\  \star  & \Sigma \rix,
\]
where $\Sigma  \in \mathbb C^{(n_1+n_3+n_4-r), (n_1+n_3+n_4-r)}$ is nonsingular. Take
\[  K=V \mat{cc} K_{11} & K_{12} \\ K_{21}^H & K_{22} \rix V^H, \quad
       F_H=V \mat{cc} 0 & 0 \\ 0 & F_{22}^{(4)}  \rix V^H, \quad F_S=0. \]
We then have that
\[
\rank(E+BKB^H)=r, \quad E+BKB^H\geq 0, \quad R+BF_HB^H\geq 0,
\]
and the pair
$(E+BKB^H, J-(R+BF_HB^H))$ is regular and of index at most one.
\eproof

The following corollary characterizes the case that the rank of $E+BKB^H$ is maximized.

\begin{corollary}\label{corollary1} Consider a pHDAE system of the form (\ref{phdae}). There exists a matrix $K$ such that
\[
E+BKB^H\geq 0,
\]
the pair $(E+BKB^H, J-R)$ is regular and of index at most one and
 \[
 \rank(E+BKB^H)=\rank  \mat{cc} E & B \rix=\max_{\hat K \in \mathbb C^{m, m}} \rank (E+B\hat K B^H),
 \]
if and only if
\begin{equation}\label{cond11}
\rank \mat{ccc} E & (J-R)\cS_\infty\left (\mat{c} E \\ B^H \rix\right ) & B \rix=n.
\end{equation}
\end{corollary}
\proof By the sufficiency proof of Theorem \ref{theorem5}, there exists a matrix $K$ such that $E+BKB^H\geq 0$,
the pair $(E+BKB^H, J-R)$ is regular and of index at most one and
 \[
 \rank(E+BKB^H)=\rank  \mat{cc} E & B \rix=\max_{\hat K \in \mathbb C^{m, m}} \rank (E+B\hat K B^H)
 \]
if and only if  conditions (\ref{cond1}) and (\ref{cond3}) hold, and
\[  n_3+n_4=\rank \mat{cc} E & B \rix, \]
and thus, if and only if
\[ n_6=n_1+n_2=0, \]
or equivalently, condition (\ref{cond11}) holds. \eproof

We can also combine regularization, index reduction and stabilization via proportional and derivative output feedback.

\begin{theorem}\label{stability-2}  Consider a pHDAE system of the form (\ref{phdae}). There exist feedback matrices  $K$, $F_S=-F_S ^H$, $F_H=F_H ^H$ such that
the   pair $(E+BKB ^H, J+BF_SB ^H-(R+BF_HB ^H))$ is regular, of index at most one, has all its finite eigenvalues in the open left {\color{black} complex half plane}, and \begin{equation}\label{eqA1}
    E+BKB^H\geq 0, \quad R+BF_HB ^H\geq 0
\end{equation}
if and only if conditions (\ref{cond1}), (\ref{cond3}), and (\ref{con-S1}) for all purely imaginary $s$, hold.
Moreover, under these conditions
for a given integer $r$, there exist matrices $K$,
$F_S=-F_S ^H$ and $F_H=F_H ^H$ such that
the pair $(E+BKB ^H, (J+BF_SB ^H)-(R+BF_HB ^H))$ is regular, of index at most one, has all its finite eigenvalues in the open left complex half plane, (\ref{eqA1}) holds, and
\[ \rank(E+BKB ^H)=r
\]
if and only if   (\ref{cond9}) holds.
\end{theorem}

\proof  The necessity of conditions (\ref{cond1}), (\ref{cond3})  and (\ref{cond9}) follow from Theorem \ref{theorem6} and the condition (\ref{con-S1}) is a standard condition in linear control~\cite{Kai80}.

For the sufficiency, for any integer $r$ satisfying (\ref{cond9}), let $K=K^H$ and $F_H\geq 0$  be chosen  as in the sufficiency proof of Theorem \ref{theorem6}, i.e., such that
the pair $(E+BKB ^H, J-(R+BF_HB ^H))$ is regular and of index at most one,
\[
E+BKB^H\geq 0, \quad R+BF_HB^H\geq 0, \quad \rank(E+BKB^H)=r.
\]
Let $\tilde F_H\geq 0$ be such that
\[
\rank (R+B(F_H+\tilde F_H)B ^H)=\rank \mat{cc} R+BF_HB ^H & B \rix=\rank \mat{cc} R & B \rix.
\]
Note that for all purely imaginary $s$ we have that
\[
\rank \mat{cc} J-(R+BF_HB ^H)-s(E+BKB ^H) & B \rix = \rank \mat{cc} J-R-sE & B\rix=n,
\]
and it follows from the sufficiency proof of Theorem \ref{stability-1} that the  pair
\[
(E+BKB ^H, J-B(F_H+\tilde F_H)B ^H)
\]
has all its finite eigenvalues in the open left {\color{black} complex half plane}. Furthermore,
\[ R+B(F_H+\tilde F_H)B^H=R+BF_HB^H+B\tilde F_HB^H\geq 0,
\]
and
\[ \cT_\infty(E+BKB ^H)=\cS_\infty(E+BKB ^H),
\]
and by Lemma \ref{lemma-S1} it follows that
\begin{eqnarray*}
& &  \cT_\infty ^H(E+BKB ^H) (J-(R+B(F_H+\tilde F_H)B^H))\cS_\infty(E+BKB^H)  \\
&& \quad  =\cT_\infty^H(E+BKB^H) (J-(R+BF_HB^H))\cS_\infty(E+BKB ^H) \\
 & & \quad \quad  -\cT_\infty ^H(E+BKB^H) (B\tilde F_HB^H)\cS_\infty(E+BKB^H)
\end{eqnarray*}
is nonsingular. Therefore, the  pair  $(E+BKB ^H, J-B(F_H+\tilde F_H)B ^H)$ is of index at most one.
\eproof

\begin{remark} \label{rem:derfeed}{\rm
Consider the condensed form (\ref{condensed-form2}). Then for $K=V\mat{cc} K_{11} & K_{12} \\ K_{12}^H & K_{22} \rix V^H$ the closed loop system
$(E+BKB^H, J-R)$ has all its finite eigenvalues in the open left {\color{black} complex half plane} if and only if the  pair
\[
\left (\mat{cc} E_{33}+B_{32}K_{22}B_{32}^H & 0 \\ 0 & E_{44} \rix, \mat{cc} A_{33} & A_{34} \\ A_{43} & A_{44} \rix \right )
\]
has all its finite eigenvalues in the open left {\color{black} complex half plane}. So,  the stabilization of the pHDAE system  (\ref{phdae})  by only derivative output feedback cannot be achieved in general.
}
\end{remark}
\section{Concluding Remarks}\label{sec:conclusions}
In this paper, new characterizations have been derived for the regularization, index reduction and stabilization of port-Hamiltonian descriptor systems (\ref{phdae}) {\color{black} by proportional and derivative
output feedback} while preserving the port-Hamiltonian structure. 
%
Future work will include the development and implementation of numerical methods for optimal robust output feedback stabilization.

\section*{Appendix}

\textbf{Constructive proof of Lemma~\ref{lemma1}.}
In this proof, we use QR decompositions and singular value decompositions, see \cite{GolV96} to determine the mentioned unitary matrices.

Step 1. Determine a unitary matrix $U_1$ such that
\[ U_1^HB=\bmat{ &          \cr
              \mu_1 & B_1    \cr
           n-\mu_1 & 0       \cr}, \]
where $\rank(B_1)=\mu_1$ and set
\[ U_1^HEU_1=\bmat{ & \mu_1 & n-\mu_1 \cr
                     \mu_1   & E_{11}^{(1)}        & E_{12}^{(1)} \cr
                 n-\mu_1    & (E_{12}^{(1)})^H & E_{22}^{(1)} \cr}, \]
and
\[  U_1^H (J-R) U_1=\bmat{ & \mu_1 & n-\mu_1 \cr
    \mu_1   & J_{11}^{(1)}-R_{11}^{(1)}              & J_{12}^{(1)}-R_{12}^{(1)} \cr
    n-\mu_1    & -(J_{12}^{(1)})^H-(R_{12}^{(1)})^H & J_{22}^{(1)}-R_{22}^{(1)}  \cr},
\]
where $E_{22}^{(1)}\geq 0$ since $E\geq 0$.

Step 2. Determine a unitary matrix $U_2$  such that
\[
U_2^H E_{22}^{(1)} U_2 =\bmat{ & \mu_2 & n-\mu_1-\mu_2 \cr
                        \mu_2    & \hat E_{22} & 0 \cr
                n-\mu_1-\mu_2  & 0                & 0  \cr},
\]
where $\hat E_{22}>0$ and set
\[ U_2^H (J_{22}^{(1)}-R_{22}^{(1)}) U_2=
\bmat{ & \mu_2 & n-\mu_1-\mu_2 \cr
  \mu_2       & J_{22}^{(2)}-R_{22}^{(2)}               & J_{23}^{(2)}-R_{23}^{(2)} \cr
   n-\mu_1-\mu_2    & -(J_{23}^{(2)})^H-(R_{23}^{(2)})^H & J_{33}^{(2)}-R_{33}^{(2)} \cr}.
\]

Step 3. Determine a unitary matrix $U_3$ such that
\[
U_3^H ( J_{33}^{(2)}-R_{33}^{(2)})=\bmat{ &                       \cr
        n_5   & \tilde J_3-\tilde R_3     \cr
        n_6   & 0                    \cr},
\]
where $\rank(\tilde J_3-\tilde R_3)=n_5$. Set
\[
U_3^H J_{33}^{(2)} U_3=\bmat{ & n_5 & n_6 \cr
         n_5  & J_{55} & J_{56}^{(3)} \cr
        n_6  & -(J_{56}^{(3)})^H & J_{66}^{(3)} \cr}, \quad
   U_3 R_{33}^{(2)} U_3^H=\bmat{ & n_5 & n_6 \cr
    n_5  & R_{55} & R_{56}^{(3)} \cr
   n_6  & (R_{56}^{(3)})^H & R_{66}^{(3)} \cr}.
\]
Then
\begin{eqnarray*}
U_3^H( J_{33}^{(2)}-R_{33}^{(2)}) U_3&=&\bmat{ & n_5 & n_6 \cr
    n_5  & J_{55}-R_{55} & J_{56}^{(3)}-R_{56}^{(3)} \cr
    n_6  & -(J_{56}^{(3)})^H-(R_{56}^{(3)})^H  & J_{66}^{(3)}-R_{66}^{(3)}   \cr}\\
    &=&
    \bmat{ & n_5 & n_6 \cr
   n_5  & J_{55}-R_{55} & J_{56}^{(3)}-R_{56}^{(3)} \cr
  n_6  & 0                    & 0                                       \cr}.
  \end{eqnarray*}
Note that $R\geq 0$ and $J=-J^H$, so, $R_{33}^{(2)}\geq 0$, $J_{33}^{(2)}=-(J_{33}^{(2)})^H$, and thus,
\[
R_{66}^{(3)} = J_{66}^{(3)}=0, \quad R_{56}^{(3)}=0, \quad J_{56}^{(3)}=0.
\]
 Define
\[
\tilde U_1=\mat{cc} I & 0\\ 0& U_3 \rix \mat{cc} I & 0\\ 0& U_2 \rix U_1,
\]
then
%
\[
\tilde U_1^H B=\bmat{&      \cr
               \mu_1 & B_1 \cr
               \mu_2 & 0     \cr
                   n_5 & 0     \cr
                   n_6 & 0    \cr}, \quad
   \tilde U_1^H E\tilde U_1  =\bmat{ & \mu_1  & \mu_2         & n_5 & n_6 \cr
             \mu_1 & \hat E_{11}     & \hat E_{12} &  0   & 0   \cr
             \mu_2 & \hat E_{12}^H & \hat E_{22} &  0   & 0   \cr
               n_5   & 0                    & 0                & 0    & 0    \cr
               n_6   & 0                    & 0                & 0    & 0     \cr}, \]
\[ \tilde U_1^H(J-R) \tilde U_1=\bmat{ & \mu_1                                         & \mu_2                                     & n_5                             & n_6 \cr
\mu_1 & \hat J_{11}-R_{11}                  & \hat J_{12}-\hat R_{12}            &  \hat J_{13}-\hat R_{13}   & \hat J_{14} -\hat R_{14}  \cr
\mu_2 & -\hat J_{12}^H-\hat R_{12}^H & \hat J_{22}-\hat R_{22}             &  \hat J_{23}-\hat R_{23}   & \hat J_{24}-\hat R_{24}    \cr
n_5    & -\hat J_{13}^H-\hat R_{13}^H & -\hat J_{23}^H-\hat R_{23}^H   & J_{55}-R_{55}                  & 0                \cr
n_6    & -\hat J_{14}^H-\hat R_{14}^H     & -\hat J_{24}^H-\hat R_{24}^H      & 0                                     & 0     \cr},
\]
where
\[
\rank(B_1)=\mu_1, \ \rank(J_{55}-R_{55})=n_5, \  \rank(\hat E_{22})=\mu_2.
\]
In addition, using $R\geq 0$, we also have that
\[
\hat R_{14}=0, \quad \hat R_{24}=0,
\]

Step 4. Construct unitary  matrices $U_4$ and $V$ such that
\begin{eqnarray*}
& & U_4^H \mat{cc} \hat E_{11} & \hat E_{12} \\ \hat E_{12}^H & \hat E_{22} \rix U_4=
     \bmat{ & n_1 & n_2 & n_3 & n_4 \cr
         n_1 & E_{11}     & E_{12}     & E_{13}      & E_{14} \cr
         n_2 & E_{12}^H & E_{22}     & E_{23}      & E_{24} \cr
         n_3 & E_{13}^H & E_{23}^H & E_{33}     & E_{34} \cr
         n_4 & E_{14}^H & E_{24}^H & E_{34}^H & E_{44} \cr}, \\
& & U_4^H \mat{c} B_1 \\ 0 \rix V =\bmat{ & m-n_3 & n_3 \cr
    n_1   & 0         & B_{12} \cr
    n_2   & B_{21} & B_{22} \cr
     n_3   & 0         & B_{32} \cr
    n_4   & 0         & 0         \cr}, \quad
      U_4^H \mat{c} \hat J_{14} \\ \hat J_{24} \rix=
      \bmat{ & n_6     \cr
   n_1 & J_{16} \cr
   n_2 & J_{26} \cr
   n_3 & 0         \cr
   n_4 & 0         \cr},
\end{eqnarray*}
where
\[
\rank \mat{c} J_{16} \\ J_{26} \rix=n_1+n_2, \ \rank(B_{21})=n_2, \ \rank(B_{32})=n_3,
\]
and
\[
\rank \mat{cccc} E_{14}^H & E_{24}^H & E_{34}^H & E_{44} \rix=n_4,
\]
which, together with $E\geq 0$, yields $ E_{44}>0$. Moreover,
\begin{eqnarray*} \rank \mat{cccccc} E_{11}     & E_{12}     & E_{13}      & E_{14} & 0         & B_{12} \\
 E_{12}^H & E_{22}     & E_{23}      & E_{24} & B_{21} & B_{22} \\
                              E_{13}^H & E_{23}^H & E_{33}      & E_{34} & 0         & B_{32} \\
                              E_{14}^H & E_{24}^H & E_{34}^H & E_{44} & 0         & 0 \rix
  &=&\rank \mat{ccc} \hat E_{11} & \hat E_{12} & B_1 \\ \hat E_{12}^H & \hat E_{22} & 0 \rix\\ &=&\mu_1+\mu_2=n_1+n_2+n_3+n_4. \end{eqnarray*}
Then
\[
U=\mat{cc} U_4 & 0\\ 0& I \rix \tilde U_1,
\]
and $V$ are the transformation matrices to the form (\ref{condensed-form1}).

\end{document}